\documentclass[11pt,oneside,english]{amsart}
\usepackage[T1]{fontenc}
\usepackage[latin9]{inputenc}
\usepackage[letterpaper]{geometry}
\usepackage{amsthm}
\usepackage{setspace}
\usepackage{amssymb}
\makeatletter
\numberwithin{equation}{section}
\numberwithin{figure}{section}
\theoremstyle{plain}
\newtheorem{thm}{Theorem}
  \theoremstyle{plain}
  \newtheorem{lem}[thm]{Lemma}
  \theoremstyle{plain}
  \newtheorem{cor}[thm]{Corollary}
  \theoremstyle{remark}
  \newtheorem{rem}[thm]{Remark}
  \theoremstyle{definition}
  \newtheorem*{example*}{Example}

\makeatother

\usepackage{babel}

\begin{document}

\title{Growth properties of the Fourier transform}

\author{William O. Bray}

\author{Mark A. Pinsky}

\address{Department of Mathematics \& Statistics, University of Maine, Orono,
Maine 04469}

\email{bray@math.umaine.edu}

\address{Department of Mathematics, Northwestern University, Evanston, Illinois
60208}

\email{mpinsky@math.northwestern.edu}

\dedicatory{In memory of Professor \"Caslav V. Stanojevi\'c as teacher, mentor,
and friend.}
\begin{abstract}
In a recent paper by the authors, growth properties of the Fourier
transform on Euclidean space and the Helgason Fourier transform on
rank one symmetric spaces of non-compact type were proved and expressed
in terms of of a modulus of continuity based on spherical means. The
methodology employed first proved the result on Euclidean space and
then, via a comparison estimate for spherical functions on rank one
symmetric spaces to those on Euclidean space, we obtained the results
on symmetric spaces. In this note, an analytically simple, yet overlooked
refinement of our estimates for spherical Bessel functions is presented
which provides significant improvement in the growth property estimates. 
\end{abstract}
\maketitle

\section{Euclidean Space}

In a paper by the authors \cite{BP}, we proved the following result
providing a growth property of the Euclidean Fourier transform.
\begin{thm}
\label{thm:bp-08}Let $1\leq p\leq2$ and $n\geq2$. Then there is
a constant $C=C(p,n)$ such that the following hold.~\end{thm}
\begin{enumerate}
\item If $f\in L^{p}(\mathbb{R}^{n})$ and $1<p\leq2$, then\[
\left(\int_{\mathbb{R}^{n}}\min\left\{ 1,\left(\frac{|\xi|}{r}\right)^{2q}\right\} |\widehat{f}(\xi)|^{q}d\xi\right)^{1/q}\leq C\,\Omega_{p}[f](\frac{1}{r}).\]

\item If $f\in L^{1}(\mathbb{R}^{n}),$ then\[
\sup_{\xi}\left[\min\left\{ 1,\left(\frac{|\xi|}{r}\right)^{2}\right\} |\widehat{f}(\xi)|\right]\leq C\,\Omega_{1}[f](\frac{1}{r}).\]

\end{enumerate}
\noindent Here the modulus of continuity was defined as\[
\Omega_{p}[f](r)=\sup_{0\leq t\leq r}\Vert M^{t}f-f\Vert_{p},\]
and $M^{t}f$ is the usual spherical mean of $f$,\[
M^{t}f(x)=\frac{1}{\omega_{n-1}}\int_{S^{n-1}}f(x+t\omega)\, d\omega,\]
$S^{n-1}$ is the unit sphere in $\mathbb{R}^{n}$, $\omega_{n-1}$
its total surface measure with respect to the usual induced measure
$d\omega$. In essence, the proof of this theorem is based on three
things: (1) the Fourier transform identity $\widehat{M^{t}f}(\xi)=\widehat{f}(\xi)\, j_{\frac{n-2}{2}}(t|\xi|)$,
where $j_{\alpha}(r)$ is normalized spherical Bessel function of
order $\alpha$,\[
j_{\alpha}(r)=2^{\alpha}\Gamma(\alpha+1)r^{-\alpha}J_{\alpha}(r),\]
(2) the Hausdorff-Young theorem, and (3) a careful estimate of (see
\cite{BP}, Lemma 6) $I(\lambda/r,z)=1-j_{\alpha}(\lambda z/r)$ of
the form,\begin{equation}
C_{1,\alpha}\min\left\{ 1,\left(\tfrac{\lambda}{r}\right)^{2}\right\} \leq\int_{0}^{1}I\left(\tfrac{\lambda}{r},z\right)dz\leq\sup_{0\leq z\leq1}I\left(\tfrac{\lambda}{r},z\right)\leq C_{2,\alpha}\min\left\{ 1,\left(\tfrac{\lambda}{r}\right)^{2}\right\} .\label{eq:old-est}\end{equation}
Here, $C_{k,\alpha}$ are positive constants.

A technically simple refinement of this estimate leads to the following
generalization; a result of the same form as Theorem \ref{thm:bp-08},
yet dispenses with the need for the supremum in the spherical modulus
of continuity.
\begin{thm}
\label{thm:bp-09}Let $1\leq p\leq2$ and $n\geq2$. Then there exists
a constant $C=C(p,n)$ such that the following hold.~\end{thm}
\begin{enumerate}
\item If $1<p<2$ and $f\in L^{p}(\mathbb{R}^{n}),$ then\[
\left(\int_{\mathbb{R}^{n}}\min\left\{ 1,(t|\xi|)^{2q}\right\} |\widehat{f}(\xi)|^{q}d\xi\right)^{1/q}\leq C\,\Vert M^{t}f(\cdot)-f(\cdot)\Vert_{p}.\]

\item If $f\in L^{1}(\mathbb{R}^{n})$, then\[
\sup_{\xi}\left[\min\left\{ 1,(t|\xi|)^{2}\right\} \right]|\widehat{f}(\xi)|\leq C\,\Vert M^{t}f(\cdot)-f(\cdot)\Vert_{1}.\]

\item If $f\in L^{2}(\mathbb{R}^{n})$, then the result takes sharper form:\[
\left(\int_{\mathbb{R}^{n}}\min\{1,(t|\xi|)^{4}\}|\widehat{f}(\xi)|^{2}d\xi\right)^{1/2}\asymp\Vert M^{t}f(\cdot)-f(\cdot)\Vert_{2}.\]
($r(t)\asymp s(t)$ means the left hand side is bounded above and
below by positive constants times the right hand side)
\end{enumerate}
The sharper form in the the case $p=2$ is because the inequality
in the Hausdorff-Young theorem for $p<2$ becomes equality in the
Plancherel theorem. The proof of this result follows the same method
as for Theorem \ref{thm:bp-08}, with the estimate for spherical Bessel
functions given above (Lemma 6 in \cite{BP}) replaced by the following.
\begin{lem}
Let $\alpha>-\tfrac{1}{2}$. Then there are positive constants $c_{1,\alpha}$
and $c_{2,\alpha}$ such that\begin{equation}
c_{1,\alpha}\min\{1,(\lambda t)^{2}\}\leq1-j_{\alpha}(\lambda t)\leq c_{2,\alpha}\min\{1,(\lambda t)^{2}\}\label{eq:basic est}\end{equation}
for all $\lambda>0$.\end{lem}
\begin{proof}
The proof makes use the the Mehler formula for the spherical Bessel
function given by\[
j_{\alpha}(\lambda t)=\frac{2\Gamma(\alpha+1)}{\sqrt{\pi}\Gamma(\alpha+\frac{1}{2})}\int_{0}^{1}(1-y^{2})^{\alpha-\frac{1}{2}}\cos(\lambda ty)\, dy.\]
It follows that\[
1-j_{\alpha}(\lambda t)=\frac{4\Gamma(\alpha+1)}{\sqrt{\pi}\Gamma(\alpha+\frac{1}{2})}\int_{0}^{1}(1-y^{2})^{\alpha-\frac{1}{2}}\sin^{2}\left(\frac{\lambda ty}{2}\right)dy.\]
Since, $\sin\frac{\lambda ty}{2}\geq\frac{\lambda ty}{\pi}$, provided
$\lambda t\leq\pi$, it follows that\[
1-j_{\alpha}(\lambda t)\geq\frac{4\Gamma(\alpha+1)}{\pi^{3/2}\Gamma(\alpha+\frac{1}{2})}(\lambda t)^{2}\int_{0}^{1}(1-y^{2})^{\alpha-\frac{1}{2}}y^{2}dy=\frac{(\lambda t)^{2}}{\pi(\alpha+1)},\]
the last step by evaluating the beta integral and simplifying the
resulting gamma functions. For all $\lambda t>0$, $|j_{\alpha}(\lambda t)|<1$.
Hence, for $\lambda t\geq\pi$, there is a constant $c>0$ such that
$1-j_{\alpha}(\lambda t)\geq c$. Combining the estimates gives the
left hand side of \eqref{eq:basic est}. The right hand side follows
by similar technique.
\end{proof}
The following corollary represents a quantified Riemann-Lebesgue lemma
and is an extension/variant of results in one dimension given in Titchmarsh
\cite[page 117]{T}.
\begin{cor}
Let $1\leq p<2$ and $n\geq2$. Then there is a positive constant
$C=C(p,n)$ such that the following hold.~\end{cor}
\begin{enumerate}
\item If $1<p<2$ and $f\in L^{p}(\mathbb{R}^{n})$, then\[
\left(\int_{|\xi|>1/t}|\widehat{f}(\xi)|^{q}d\xi\right)^{1/q}\leq C\,\Vert M^{t}f(\cdot)-f(\cdot)\Vert_{p}.\]

\item If $f\in L^{1}(\mathbb{R}^{n}),$ then\[
\sup_{|\xi|>1/t}|\widehat{f}(\xi)|\leq C\,\Vert M^{t}f(\cdot)-f(\cdot)\Vert_{1}.\]

\item If $f\in L^{2}(\mathbb{R}^{n})$, then\[
\left(\int_{|\xi|>1/t}|\widehat{f}(\xi)|^{2}d\xi\right)^{1/2}\approx\Vert M^{t}f(\cdot)-f(\cdot)\Vert_{2}.\]
\end{enumerate}
\begin{rem}
Theorem \eqref{thm:bp-09} was also obtained by Ditzian \cite{D}
as a consequence of a rather technical result in approximation theory.
Our proof lies completely within the framework of harmonic analysis
and lends its self to the extensions described below.
\end{rem}

\section{Rank One Symmetric Spaces}

In this section we follow the notation given in \cite{BP}; basic
references for the background material are Helgason's books \cite{H1,H2}
and Koornwinder's survey paper on Jacobi functions \cite{K}. Let
$X=G/K$ where $G$ is a connected non-compact semisimple Lie group
with finite center and real rank one and $K$ is a maximal compact
subgroup. At the Lie algebra level, the Cartan decomposition has form
$\mathfrak{g}=\mathfrak{k}+\mathfrak{p}$, where $\mathfrak{k}$ is
the Lie algebra of $K$. The Iwasawa decomposition takes the form
$\mathfrak{g}=\mathfrak{k}+\mathfrak{a}+\mathfrak{n}$, where $\mathfrak{a}$
is a maximal abelian subalgebra of $\mathfrak{p}$ and $\mathfrak{n}$
is a nilpotent subalgebra of $\mathfrak{g}$. The rank one condition
is that $\dim\mathfrak{a}=1$. The nilpotent subalgebra $\mathfrak{n}$
has root space decomposition $\mathfrak{n}=\mathfrak{n}_{\gamma}+\mathfrak{n}_{2\gamma}$,
where $\gamma$ and $2\gamma$ are the positive roots. Let $m_{\gamma}$
and $m_{2\gamma}$ be the respective root space dimensions and set
$\rho=\frac{1}{2}(m_{\gamma}+2m_{2\gamma})$. Choose $H_{0}\in\mathfrak{a}$
such that $\gamma(H_{0})=1$. This allows identifying $\mathfrak{a}$
with $\mathbb{R}$ by the map $\mathbb{R}\ni t\rightarrow tH_{0}\in\mathfrak{a}$,
and on the dual side, $\mathfrak{a}_{\mathbb{C}}^{*}$ with $\mathbb{C}$.
At the group level, the Iwasawa decomposition has form $G=KAN$, and
we write $G\ni g=k\exp(H(g))\, n$, where $H(g)\in\mathfrak{a}$ and
$\exp$ is the exponential function. Because of the above identification,
we often write $a_{t}=\exp(H(g))$, $t\in\mathbb{R}$ being identified
with $H(g)$.

The polar decomposition of $G$ takes the form $G=KA^{+}K$, where
$A^{+}=\{a_{t}\,:\, t\geq0\}$. Following standard practice, functions
$f$ on $X$ are identified with right $K-$invariant functions on
$G$ and write $f(x)=f(g)$, where $x=gK$. In terms of this decomposition,
the invariant measure $dx$ on $X$ has the form\[
dx=\Delta(t)\, dtdk,\]
where $\Delta(t)=\Delta_{(\alpha,\beta)}(t)=(2\sinh t)^{2\alpha+1}(2\cosh t)^{2\beta+1}$,
$\alpha=(m_{\gamma}+m_{2\gamma}-1)/2$ and $\beta=(m_{2\gamma}-1)/2$,
and $dk$ is normalized Haar measure on $K$. The Laplacian on $X$
is denoted $L$ and its radial part is given by\[
L_{r}=\frac{d^{2}}{dt^{2}}+\frac{\Delta'(t)}{\Delta(t)}\frac{d}{dt}.\]
The spherical function on $X$ is the unique radial solution to the
equation\[
Lu=-(\lambda^{2}+\rho^{2})u\]
which is one at the origin of $X$. Let $M$ be the centralizer of
$A$ in $K$ and set $B=K/M$. For $x=gK\in X$ and $b=kM\in B$,
let $A(x,b)=-H(g^{-1}k)$ (called the horocycle distance function).
Then the Harish-Chandra formula for the spherical function is\[
\phi_{\lambda}(x)=\int_{B}e^{(i\lambda+\rho)A(x,b)}db,\]
where $db$ is normalized measure on $B$. If we write $x=ka_{t}K$,
then it is well known that $\phi_{\lambda}(x)=\phi_{\lambda}^{(\alpha,\beta)}(t)$,
where $\phi_{\lambda}^{(\alpha,\beta)}(t)$ is Jacobi function of
the first kind. Key properties of Jacobi functions are given in the
following three bullet items \cite{K}.
\begin{itemize}
\item $|\phi_{\mu+i\eta}^{(\alpha,\beta)}(t)|\leq\phi_{i\eta}^{(\alpha,\beta)}(t)\leq1$
for $\mu\in\mathbb{R}$ and $|\eta|\leq\rho$.
\item $|\phi_{\mu+i\eta}^{(\alpha,\beta)}(t)|\leq e^{|\eta|t}\phi_{0}^{(\alpha,\beta)}(t)\leq C(1+t)e^{(|\eta|-\rho)t}$.
\item Let $1<p<2$, and define $D_{p}=\{\lambda=\mu+i\eta\,:\,|\eta|<(\frac{2}{p}-1)\rho\}.$
Then\[
\lambda\in D_{p}\implies\phi_{\lambda}^{(\alpha,\beta)}\in L^{q}(\mathbb{R}^{+},\Delta_{(\alpha,\beta)}(t)dt),\]
where $q$ is the H\"older conjugate index: $\frac{1}{p}+\frac{1}{q}=1$.
\end{itemize}
In sharp contrast to the Euclidean space setting, the third property
written out for symmetric space states: for $\lambda\in D_{p}$, the
spherical function $\phi_{\lambda}\in L^{q}(X)$.

In \cite{BP}, we proved the following result.
\begin{lem}
\label{lem:comparison}Let $\alpha>-1/2$, $-1/2\leq\beta\leq\alpha$,
and let $t_{0}>0$. Then for $|\eta|\leq\rho$, there exists a positive
constant $C=C(\alpha,\beta,t_{0})$ such that\[
|1-\phi_{\mu+i\eta}^{(\alpha,\beta)}(t)|\geq C\,[1-j_{\alpha}(\mu t)],\]
for all $0\leq t\leq t_{0}$.
\end{lem}
In the symmetric space realm, the above gives a local estimate involving
the spherical function on $X$ with that for the spherical function
on the Euclidean tangent space to $X$ at the origin and is the technical
heart of the extension of Theorem 1 to symmetric spaces. The following
example illustrates the essential ideas underlying this estimate.
\begin{example*}
Consider the case $\alpha=\frac{1}{2}$, $\beta=-\frac{1}{2}$. Then
the Jacobi function is elementary\[
\phi_{\lambda}(t)=\phi_{\lambda}^{(1/2,-1/2)}(t)=\frac{\sin\lambda t}{\lambda\sinh t},\]
and gives the spherical function on three dimensional real hyperbolic
space $SO_{e}(1,3)/SO(3)$. Using the fundamental theorem of calculus\[
1-\phi_{\lambda}(t)=\frac{1}{\sinh t}\int_{0}^{t}[\cosh s-\cos\lambda s]ds.\]
Substituting $\lambda=\mu+i\eta$, applying the addition theorem for
the cosine function, and the fact that modulus dominates the real
part, we obtain\begin{align*}
\left|1-\phi_{\lambda}(t)\right| & \geq\frac{1}{\sinh t}\int_{0}^{t}[\cosh s-\cosh\eta s\cos\mu s]ds.\end{align*}
It follows that\begin{eqnarray*}
\left|1-\phi_{\lambda}(t)\right| & \geq & \frac{1}{\sinh t}\int_{0}^{t}\left[1-\frac{\cosh\eta s}{\cosh s}\cos\mu s\right]ds\\
 & \geq & \frac{t}{\sinh t}\frac{1}{t}\int_{0}^{t}[1-\cos\mu s]ds\\
 & = & \frac{t}{\sinh t}\left[1-j_{1/2}(\mu t)\right].\end{eqnarray*}
For fixed $t_{0}>0$, the ratio $\frac{t}{\sinh t}\geq C$ for all
$0\leq t\leq t_{0}$, which gives the proof of the lemma for this
example. The proof in the general case is based on more elaborate
estimates applied to the Mehler identity for the Jacobi functions
(the Mehler identity for the example is straightforward via the fundamental
theorem of calculus).
\end{example*}
As a consequence of the second bullet item above, \[
\lim_{t\rightarrow\infty}\phi_{\lambda}^{(\alpha,\beta)}(t)=0\]
uniformly on any strip of the form $\{\lambda=\mu+i\eta\,:\,\mu\in\mathbb{R},\,|\eta|\leq\eta_{0}<\rho\}$.
Combining the above lemma and this fact with \eqref{eq:basic est}
gives the following estimate.
\begin{lem}
\label{lem:sym space impr est}Let $\alpha>-1/2$, $-1/2\leq\beta\leq\alpha$,
and let $0<\eta_{0}<\rho$. Then there exists a positive constant
$C=C(\alpha,\beta,\eta_{0})$ such that\[
|1-\phi_{\mu+i\eta}^{(\alpha,\beta)}(t)|\geq C\,\min\{1,(\mu t)^{2}\}\]
for all $\mu\in\mathbb{R}$ , $|\eta|\leq\eta_{0}$, and $t>0$.
\end{lem}
The Helgason Fourier transform for functions defined on $X$ is given
by\[
\widehat{f}(\lambda,b)=\int_{X}f(x)\, e^{(-i\lambda+\rho)A(x,b)}dx.\]
The following estimate due to Sarkar and Sitaram \cite{SS}for this
transform provides one avenue for developing an analog of Theorem
\ref{thm:bp-09} in symmetric spaces; it is a direct consequence of
the aforementioned integrability property of Jacobi functions.
\begin{lem}
\label{lem:sarkar}Let $1\leq p<2$ and $f\in L^{p}(X)$. Then for
$\lambda=\mu+i\eta\in D_{p}$, $\widehat{f}(\lambda,b)$ is defined
a.e. ($b$) and\begin{equation}
\int_{B}|\widehat{f}(\mu+i\eta,b)|\, db\leq c_{p}(|\eta|)\,\Vert f\Vert_{p},\label{eq:sakar}\end{equation}
where $c_{p}(\cdot)$ is a positive function defined on $[0,(\frac{2}{p}-1)\rho)$.
\end{lem}
Group theoretically, the spherical mean of a function $f$ on $X$
is given by\[
M^{t}f(g)=\int_{K}f(gka_{t})\, dk.\]
The main result generalizing Theorem 12 of \cite{BP} is the following.
\begin{thm}
\label{thm:Lp-result}Let $1\leq p<2$ and $f\in L^{p}(X)$. Then
for $|\eta|<(2/p-1)$ and $t>0$, there exists a positive function
$c_{p}(|\eta|)$ such that\[
\sup_{\mu}\left[\min\{1,(t\mu)^{2}\}\int_{B}|\widehat{f}(\mu+i\eta,b)|\, db\right]\leq c_{p}(|\eta|)\,\Vert M^{t}f(\cdot)-f(\cdot)\Vert_{p}.\]
\end{thm}
\begin{proof}
For completeness, we sketch the proof. From the operational property\[
\widehat{M^{t}f}(\lambda,b)=\phi_{\lambda}(a_{t})\widehat{f}(\lambda,b)\]
and Lemma \ref{lem:sarkar} we have\[
\left|1-\phi_{\mu+i\eta}(a_{t})\right|\int_{B}\left|\widehat{f}(\mu+i\eta,b)\right|db\leq c_{p}(|\eta|)\Vert M^{t}f(\cdot)-f(\cdot)\Vert_{p}.\]
The result then follows by applying Lemma \ref{lem:sym space impr est}. \end{proof}
\begin{cor}
Let $1\leq p<2$ and $f\in L^{p}(X)$. Then for $|\eta|<(2/p-1)\rho$
and $t>0$, there exists a positive function $c_{p}(|\eta|)$ such
that\[
\sup_{|\mu|>1/t}\int_{B}|\widehat{f}(\mu+i\eta,b)|\, db\leq c_{p}(|\eta|)\,\Vert M^{t}f(\cdot)-f(\cdot)\Vert_{p}.\]

\end{cor}
In the case $p=2$, the inequalities above breakdown. However, one
can resort to the known Plancherel theorem for the Helgason Fourier
transform and obtain a direct analog of Theorem \ref{thm:bp-09}.
The following result generalizes Theorem 14 of \cite{BP}.
\begin{thm}
\label{thm:L2-result}Let $f\in L^{2}(X)$. Then there exists a positive
constant $C$ such that\[
\left(\int_{\mathbb{R}}\min\{1,\left(\lambda t\right)^{4}\}\int_{B}|\widehat{f}(\lambda,b)|^{2}db\,|c_{X}(\lambda)|^{-2}d\lambda\right)^{1/2}\leq C\,\Vert M^{t}f(\cdot)-f(\cdot)\Vert_{2},\]
where $c_{X}(\lambda)$ is the Harish-Chandra $c-$function for $X$.\end{thm}
\begin{rem}
At the time of writing \cite{BP}, a Hausdorff-Young inequality for
the Helgason Fourier transform was known only in the case of radial
functions. Hence the analog of Theorem 11 for the case $f\in L^{p}(X)$,
$1\leq p<2$ was left as conjecture. In the same time frame as \cite{BP},
Ray and Sarkar \cite{RA} proved the Hausdorff-Young theorem and provided
nice generalizations to the inequality \eqref{eq:sakar} in the context
of Lorentz spaces using complex interpolation techniques. This has
been applied in \cite{RA2} to obtain corresponding extensions/refinements
of Theorems \ref{thm:Lp-result} and \ref{thm:L2-result} in the context
of harmonic $NA-$groups. The latter include the rank one symmetric
spaces as special cases.\end{rem}

\end{document}